\documentclass[12pt,reqno]{amsart}%
\usepackage{amsfonts,anysize}
\usepackage[charter]{mathdesign}
\usepackage[colorlinks,urlcolor=blue,citecolor=blue,linkcolor=blue]{hyperref}
\usepackage{graphicx}
\marginsize{3cm}{3cm}{1.5cm}{1.5cm}

\newtheorem{thm}{Theorem}[section]

\newtheorem{lem}[thm]{Lemma}

\theoremstyle{definition}

\theoremstyle{remark}

\numberwithin{equation}{section}
\allowdisplaybreaks
\begin{document}

\title{On the existence of vectors dual to a set of linear functionals}
\author{Shibo Liu}
\address{School of Mathematical Sciences, Xiamen University, Xiamen 361005, China}
\email{liusb@xmu.edu.cn}
\thanks{This work was supported by NSFC (11671331). It was done while the author was
visiting the Abdus Salam International Centre for Theoretical Physics. The
author is grateful to ICTP for its hospitality.}

\begin{abstract}
We give a simple proof of a crucial lemma that is established in \cite[Lemma
2.1]{MR2801386} by induction, and plays important roles in that paper and
\cite{MR2994823}.

\end{abstract}
\maketitle

In a recent paper \cite{MR2801386}, Brezhneva and Tret'yakov give an
elementary proof of the Karush--Kuhn--Tucker Theorem (that is, \cite[Theorem
1.1]{MR2801386}) in normed linear spaces. This theorem is about the
minimization problems with a finite number of inequality constraints. In their
proof of \cite[Theorem 1.1]{MR2801386}, an important step is to establish the
following lemma.

\begin{lem}
[{\cite[Lemma 2.1]{MR2801386}}]\label{l1}Let $X$ be a linear space. Let
$\xi^{i}:X\rightarrow\mathbb{R}$, $i=1,\ldots,n$, be linear functionals which
are linearly independent. Then there exists a set $\left\{  \eta_{i}\right\}
$ of $n$ linearly independent vectors in $X$ such that the matrix%
\[
A=\left(
\begin{array}
[c]{ccc}%
\langle\xi^{1},\eta_{1}\rangle & \cdots & \langle\xi^{1},\eta_{n}\rangle\\
\vdots &  & \vdots\\
\langle\xi^{n},\eta_{1}\rangle & \cdots & \langle\xi^{n},\eta_{n}\rangle
\end{array}
\right)
\]
is invertible.
\end{lem}

This lemma also plays a crucial role in \cite{MR2994823}, where an elementary
proof of the Lagrange multiplier theorem in normed linear spaces is given.

We emphasize that \ref{l1} is equivalent to the following seemingly stronger result.

\begin{thm}
\label{t1}Let $X$ be a vector space. Let $\xi^{i}:X\rightarrow\mathbb{R}$,
$i=1,\ldots,n$, be linear functionals which are linearly independent. Then
there exist $\varepsilon_{1},\ldots,\varepsilon_{n}$ in $X$ such that
$\langle\xi^{i},\varepsilon_{j}\rangle=\delta_{j}^{i}$, where $\delta_{j}^{i}$
is the Kronecker delta.
\end{thm}

In fact, let%
\[
A^{-1}=\left(
\begin{array}
[c]{ccc}%
b_{1}^{1} & \cdots & b_{n}^{1}\\
\vdots &  & \\
b_{1}^{n} & \cdots & b_{n}^{n}%
\end{array}
\right)  \text{,\qquad}\varepsilon_{j}=\sum_{k=1}^{n}b_{j}^{k}\eta_{k}%
\]
Then $AA^{-1}=I_{n}$, the $n\times n$ identity matrix, means%
\[
\langle\xi^{i},\varepsilon_{j}\rangle=\left\langle \xi^{i},\sum_{k=1}^{n}%
b_{j}^{k}\eta_{k}\right\rangle =\sum_{k=1}^{n}\langle\xi^{i},\eta_{k}\rangle
b_{j}^{k}=\left(  AA^{-1}\right)  _{j}^{i}=\delta_{j}^{i}\text{.}%
\]
Therefore, Lemma \ref{l1} and Theorem \ref{t1} are equivalent.

In \cite{MR2801386}, Lemma \ref{l1} is proved by induction on $n$. In this
note, we present a much simple proof of Theorem \ref{t1}.

\begin{proof}
[\indent Proof of Theorem \ref{t1}]We define a linear map $\varphi:X\rightarrow
\mathbb{R}^{n}$ by%
\[
\varphi(x)=\left(  \xi^{1}(x),\ldots,\xi^{n}(x)\right)  \text{.}%
\]
We claim that $\varphi$ is surjective. Otherwise, $\varphi(X)$ is a proper
subspace of $\mathbb{R}^{n}$, there exists $\lambda\in\mathbb{R}^{n}%
\backslash\left\{  0\right\}  $ which is perpendicular to $\varphi(X)$. Let
$\lambda=\left(  \lambda_{1},\ldots,\lambda_{n}\right)  $, then for all $x\in
X$,%
\[
0=\lambda\cdot\varphi(x)=\sum_{i=1}^{n}\lambda_{i}\xi^{i}(x)\text{.}%
\]
Hence the set of linear functionals $\left\{  \xi^{i}\right\}  $ is linearly
dependent, a contradiction.

Let $e_{j}=(0,\ldots,1,\ldots,0)$ be the $j$-th standard base vector of
$\mathbb{R}^{n}$ and take $\varepsilon_{j}\in\varphi^{-1}(e_{j})$,
then%
\[
\left(  \xi^{1}(\varepsilon_{j}),\ldots,\xi^{n}(\varepsilon_{j})\right)
=\varphi(\varepsilon_{j})=e_{j}\text{.}%
\]
It follows that $\xi^{i}(\varepsilon_{j})=\delta_{j}^{i}$.
\end{proof}

\end{document}